\documentclass{article}
\usepackage{indentfirst,latexsym,bm,amsmath,amsthm,amsfonts}
\begin{document}
\title{Mackey Theory for $p$-adic Lie groups}

\author{BinYong HSIE
      \\{\small{ Department of Mathematics, PeKing University, BeiJing, 100871, P.R. China
      }}
      \\{\small{E-mail: tjxb@pku.edu.cn}}}
\date{2003.4}
\maketitle{}
\begin{abstract}

This paper gives a $p$-adic analogue of the Mackey theory, which
relates representations of a group of type $G=H\times_{t} A $ to
systems of imprimitivity.

\end{abstract}


MR(2000) Subject Classification: 22E50, 20G05

\setcounter{section}{-1}

\section{\bf{Introduction}}

There are many methods for representations of the groups which are
semi-direct products. For Heisenberg groups, there are Stone-Von
Neumann Theorem and Weil's acta paper. For general cases, one can
use Jacquet functor. In this paper, we consider smooth
representations of a group $G=H\times_{t}A$ with $H$ a locally
compact and totally disconnected group and $A$ an Abelian
topological group such that $A$ and its dual $\hat{A}$ are both
locally compact and totally disconnected Abelian groups. Our
method is different from Weil's paper and doesn't need Jacquet
functor, but is along Mackey's idea. For Jacquet functor,
\cite{Bump} is a good reference.

Mackey has considered representations of Lie groups of type
$G=H\times_{t} \mathbb{R}^{n} $ with $H$ a Lie group. He relates
the representations of $G$ to systems of impremitivity Of
($H,\mathbb{R}^{n}$) (cf. Lemma 1 and Lemma 2 in our case). There
is a one to one correspondence between them. A system of
imprimitivity of $(H,\hat{A})$ means a representation $\pi$ of $H$
and a projection valued measure $P$ based on $\hat{A}$ such that
$$\pi(h)P_{E}\pi(h)^{-1}=P_{h[E]}$$
where $h \in H$ and $E$ is a Borel subset of $\widehat{A}$. In
\cite{V}, Varadarajan relates systems of impremitivity of
$(H,\hat{A})$ to ``cocycles''. The calculate of ``cocycles'' is
not an easy work. In our case, ``cocycles'' aren't needed, because
the topology is better (totally disconnected) and the
representations are smooth, i.e.,``locally constant''. We make the
sheaf theory in the sense of Bernstein and Zelevinsky instead.

Section 1 gives the spectral decomposition of a smooth
representation of $A$, which corresponds to Fourier analysis of
$\mathbb{R}^{n}$ in real case. Section 2 states Mackey Theory of
$p$-adic groups.

The sheaf theory of B-Z makes a projection valued measure
$(P;\hat{A},V)$ i.e. a $C^{\infty}_{c}(\hat{A})$-module $V$ into a
sheaf on $\hat{A}$. By it, we change the representation space into
the space of sections of the corresponding sheaf. We then change
the sections into certain ``compact'' supported functions on $H$.
Then we find that all of our representations are induced
representations.

At last, we shall point out that all of our representations are
complex and smooth, and that all functions are complex valued.

\section{Spectral decomposition}
\subsection{A fact}

Let $B$ be a compact and totally non-connected Abelian group. Let
$ C_{loc}(B) $ denote the set of all locally constant functions on
$B$, then $ \widehat{B}$, the dual of $B$, is contained in
$C_{loc}(B) $. We show that $\widehat{B}$ generates $ C_{loc}(B)
$, or more precisely that every function in $ C_{loc}(B) $ is a
finite linear sum of elements of $ \widehat{B} $. To see this , we
need only to show that $\chi_{E}$ for every open compact subset
$E$ is so. Let $B _{1} $ be an open subgroup such that $\chi_{E}$
is constant on every $xB_{1}$ for $x\in B$. We can regard
$\chi_{E}$ as a function on $B/B_{1}$ , which is a finite group.
There are $n=[B:B_{1}]$ elements $ \hat{b}_{1},...,\hat{b}_{n} $
in $ \hat{B} $ such that $ \hat{b}_{1},...,\hat{b}_{n} $ are all
characters of $B/B_{1}$. Then $\chi_{E}$ is a linear sum  of $
\hat{b}_{1},...,\hat{b}_{n} $. Therefore, we see that
$\widehat{B}$ generates $ C_{loc}(B) $.

\subsection{ Cosmooth projection valued measure and spectral decomposition}

Let $A$ and $\widehat{A}$, the dual of $A$, be locally compact and
totally non-connected Abelian groups. We will make this assumption
in the following. $A= \mathbb{Q}_{p}$ is such an example, but
$\mathbb{Q}_{p}^{\times}$ is not.

Let $U$ be a compact open subgroup of $A$, then the dual of $A/U$
is $U^{\perp}=\{\widehat{a}\in\widehat{A}|\forall u\in
U,<\widehat{a},u>=1\}$. Since $A/U$ is discrete, $U^{\perp}$ is a
compact group. Since the topology of $\widehat{A}$ is the
open-compact topology, $U^{\perp}$ is an open subgroup of
$\widehat{A}$. Therefore $U^{\perp}$ is compact open subgroup of
$\widehat{A}$.

Let $(\pi,V)$ be a smooth representation of $A$. Fix a vector
$v\in V$. There is a compact open subgroup $U$ of $A$ which fixes
$v$. Let $V^{U }$ be the subspace whose vectors are fixed by $U$.
Then $V^{U}$ is stable under $A$. We regard $A/U$ as functions on
$U^{\perp}$. By the result of section 1.1, $A/U$ generates
$C_{loc}(U^{\perp})$, so we can extend $\pi|V^{U}$ to a
representation $\pi^{U}$ of the algebra $C_{loc}(U^{\perp})$.
Write $P^{U}_{E}$ for $\pi^{U}(\chi_{E})$, where $E$ is an open
subset of $U^{\perp}$. We see that
$$P^{U}_{U^{\perp}}=\mathrm{I},P^{U}_{\phi}=0 \;\;\;\;\eqno(1.1.1) $$
$$P^{U}_{E}P^{U}_{F}=P^{U}_{E\cap F}\;\;\;\;\eqno(1.1.2)$$
$$P^{U}_{\cup E_{i}}=\sum P^{U}_{E_{i}}\;\;\;\;\eqno(1.1.3).$$
where $E,F,E_{i}$ are open subsets of $U^{\perp}$ and
$E_{i}\bigcap E_{j}=\phi $ unless $i=j$. We call such a system
$(P;U^{\perp},V^{U})$ a projection valued measure based on
$U^{\perp}$.

 We have
$$\pi(a)|V^{U}=\int_{U^{\perp}}x(a) dP^{U}(x)\;\;\;\;\eqno(1.2).$$

It is easy to see that there exists an unique
$(P;U^{\perp},V^{U})$ such that (1.2) is satisfied, by applying
the basic fact of section 1.1.

We extend $P^{U}$ to a projection valued measure based on
$\widehat{A}$ by setting $P^{U}(E)=P^{U}(E\bigcap U^{\perp})$.
Then it is easy to see that
$$P^{U}_{\widehat{A}}=\mathrm{I},P^{U}_{\phi}=0 \;\;\;\;\eqno(1.3.1) $$
$$P^{U}_{E}P^{U}_{F}=P^{U}_{E\cap F}\;\;\;\;\eqno(1.3.2)$$
$$P^{U}_{\cup E_{i}}=\sum P^{U}_{E_{i}}\;\;\;\;\eqno(1.3.3)$$
where $E,F,E_{i}$ are open subsets of $\widehat{A}$ and
$E_{i}\bigcap E_{j}=\phi $ unless $i=j$.

We can define $P^{U^{'}}$ for other compact  open subgroup $U^{'}$
of $A$ in the same way.

If $v\in V^{U}\bigcap V^{U^{'}}$, then
$$\pi(a)v=\int x(a)(dP^{U}(x))v,$$
$$\pi(a)v=\int x(a)(dP^{U^{'}}(x))v.$$

Let $\widehat{v}\in\widehat{V}$, the dual of $V$.Then
$$<\pi(a)v,\widehat{v}>=\int x(a)<dP^{U}(x)v,\widehat{v}>$$
$$=\int x(a)<dP^{U^{'}}(x)v,\widehat{v}>$$

Applying the fact in section 1.1 with $U^{\perp}\cdot U^{'\perp}$
instead of $B$ and the above formula, we see that for each $E$,
$$<P^{U}(E)v,\widehat{v}>=<P^{U}(E\cap U^{\perp}\cdot U^{'\perp})v,\widehat{v}>$$
$$=<P^{U^{'}}(E\cap U^{\perp}\cdot
U^{'\perp})v,\widehat{v}>=<P^{U^{'}}(E)v,\widehat{v}>$$

Since $\widehat{v}$ is arbitrary,
$$P^{U}(E)v=P^{U^{'}}(E)v\;\;\;\eqno(1.4).$$

By (1.4),we can patch all $P^{U}$ to a $P$ such that
$$P_{\widehat{A}}=\mathrm{I}, P_{\phi}=0\;\;\;\eqno(1.5.1)$$
$$P_{E}P_{F}=P_{E\bigcap F}\;\;\;\eqno(1.5.2)$$
$$P_{\bigcup E_{i}}=\sum P_{E_{i}}\;\;\;\eqno(1.5.3)$$
where $E,F,E_{i}$ are as in (1.3). Such a $P$, i.e. a $P$
satisfying (1.5.1), (1.5.2) and (1.5.3) is called a projection
valued measure.

Furthermore, for any given $v\in V$, there exists a compact open
subset $E(v)$ such that
$$P_{E}v=P_{E\cap E(v)}v\;\;\;\eqno(1.5.4)$$

We call a projection valued measure $(P; \hat{A},V )$ satisfying
(1.5.4) a {\bf{cosmooth projection valued measure}}.

Now, let $(P; \hat{A},V )$ be a cosmooth projection valued
measure. For any $v\in V$, define
$$\pi(a)v=\int_{E(v)}x(a)dP(x)v\;\;\;\eqno(1.6).$$

We can show that $\pi(a)v$ does not depend the choice of $E(v)$.
In fact, for another choice $E^{'}(v)$,
$$P_{E}v=P_{E\cap E(v)\cap E^{'}(v)}v$$
and therefore (1.6) become
$$\pi(a)v=\int_{E(v)\cap E^{'}(v)}x(a)dP(x)v\;\;\;\eqno(1.6.1)$$
So it is not depend on $E(v)$.

 Write (1.6) simply as
 $$\pi(a)=\int_{\widehat{A}}x(a)dP(x)\;\;\;\eqno(1.6').$$

For every compact open subgroup $E$ of $\widehat{A}$, by formula
(1.5.1)-(1.5.3), we can define
$$\pi(f)v=\int_{E}f(x)dP(x)v,\;\;\;v \in P_{E}V$$
 which defines a representation of the algebra $C_{loc}(E)$ on
 $P_{E}V$ and therefore $\pi$ is a representation of $A$ on $P_{E}V$.
 So (1.6) defines a representation of $A$ on $V$. We obtain the
 main result of this section:

$\mathbf{Theorem}$ 1. For a smooth representation $(\pi,V)$ of
$A$, there exists a unique cosmooth projection valued measure
$(P;\widehat{A},V)$ such that
$$\pi(a)v=\int_{\widehat{A}}x(a)dP(x)v\;\;\; a\in A,\; v\in V.$$
Conversely, given a cosmooth projection valued measure
$(P;\widehat{A},V)$, the above formula defines a smooth
representation $\pi$ of $A$.

Furthermore, we see that an operator on $V$ commutes with $\pi$ if
and only if it commutes with $P$.

In the next part, (1.6') is always in the sense of (1.6).

\section{Representation and system of imprimitivity}
\subsection{Semidirect product}

Let $H$ be a locally compact and totally disconnected group and
$A$ be as in section 1. Assume that there is a continuous
homomorphism $t$ of $H$ into the automorphism group of $A$. We
write $h[a]$ simply for $t_{h}(a)$. We now define a group
$G=H\times_{t}A$ by
$$(h,a)(h^{'},a^{'})=(hh^{'},at_{h}(a^{'}))\;\;\;\eqno(2.1).$$
It is easy to verify that $G$ is really a group with the identity
$e=(e_{H},e_{A})$. Furthermore
$$(h,a)^{-1}=(h^{-1},h^{-1}[a])\;\;\;\eqno(2.2).$$
$G$ is called the semidirect product of $H$ and $A$ relative to
$t$. Since $t$ is continuous,  $G$ becomes a topological group
with the product topology.

A quick calculation shows that
$$(h,a)(h^{'},a^{'})(h,a)^{-1}=(hh^{'}h^{-1},ah[a^{'}]t_{hh^{'}h^{-1}}[a^{-1}])\;\;\;\eqno(2.3).$$
It follows that $\tilde{A}=\{(e_{H},a):a\in A\}$ is a closed
normal subgroup of $G$, and that
$$(h,a)(e_{H},a^{'})(h,a)^{-1}=(e_{H},ah[a^{'}]a^{-1})\;\;\;\eqno(2.4).$$

We put
$$ \tilde{H} = \{ (h,e_{A}):h\in H \},$$
then $ \tilde{H} $ is a closed subgroup of $G$. We identity $H$
with $\tilde{H}$ and $A$ with $\tilde{A}$, then we have
$$G=AH\;\;\;\eqno(2.5.1),$$
$$\{ e  \}=A\cap H\;\;\;\eqno(2.5.2),$$
$$h[a]=hah^{-1}\;\;\;\eqno(2.5.3).$$

\subsection{Representation of $G$ and system of imprimitivity}

In this section, we relate a smooth representation of $G$ to a
cosmooth system of imprimitivity.

$\mathbf{Definition}$. Let $X$ be a continuous $H$-space. {\bf{A
system of imprimitivity}} for ($H$,$X$) acting on $V$ is a pair
$(\pi,P; V)$, where $\pi$ is a smooth representation of $H$ on $V$
and $P(E\rightarrow P_{E})$ is a projection valued measure based
on $X$, such that they satisfy a relation:
$$\pi_{h}P_{E}\pi_{h^{-1}}=P_{h.E}\;\;\;\eqno(2.6),$$
where, $h\in H$ and $E$ is an open subset of $X$. Furthermore, if
$P$ is cosmooth, then $(\pi,P;V)$ is called {\bf{cosmooth system
of imprimitivity}}.

Two systems $(\pi,P;V)$ and $(\pi^{'},P^{'};V^{'})$ based on the
same $H$-space $X$ are said to be {\bf{equivalent}} if and only if
there exists an isomorphism $T$ from $V$ to $V^{'}$ such that
$$\pi^{'}(h)=T\pi(h)T^{-1}\;\;\eqno(2.7.1),$$
$$P^{'}_{E}=TP_{E}T^{-1}\;\;\;\eqno(2.7.2),$$
where, $h\in H$ and $E$ is an open subset of $X$. We say that a
cosmooth system of imprimitivity $(\pi,P;V)$ is {\bf{irreducible}}
if and only if there is no subspace other than 0 and $V$ which is
invariant under all $P_{E}$ and $\pi_{h}$.

We define a homomorphism $t^{'}$ of $H$ to the the automorphism
group of $\widehat{A}$ by
$$t^{'}_{h}( \hat{a} )(a)=\widehat{a}(t_{h^{-1}}(a))\;\;\;\eqno(2.8),$$
and we write simply $h[\widehat{a}]$ for $t^{'}_{h}(\widehat{a})$.
Then $\widehat{A}$ becomes a continuous $H$-space.

If $\pi$ is a smooth representation of $G$, then $\pi$
restrictions to $A$ and $H$ are also smooth.

{\bf Lemma 1.} Let $\pi_{1}$ and $\pi_{2}$ be smooth
representations of $A$ and $H$ respectively in a vector space $V$,
and let $P$ be the corresponding cosmooth projection valued
measure on $\widehat{A}$ for $\pi_{1}$. Then a necessary and
sufficient condition such that there exists a smooth
representation $\pi$ of $G$ in $V$
 whose restrictions to $A$ and $H$ are $\pi_{1}$ and $\pi_{2}$
 respectively, is that $(\pi_{2},P;V)$ is a cosmooth system of
 imprimitivity for $H$ based on $\widehat{A}$. In this case, $\pi$
 is unique.

 \begin{proof}
 Let $\pi$ be a smooth representation of $G$ in $V$, and let
 $\pi_{1}$, $\pi_{2}$ be the restrictions to $A$,$H$ respectively.
 Now
 $$hah^{-1}=h[a]\;\;\;\eqno(2.9)$$
 so that
 $$\pi_{2}(h)\pi_{1}(a)\pi_{2}(h)^{-1}=\pi_{1}(h[a])\eqno(2.10)$$
 for all $(h,a)\in H \times A$. Let $P$ be the corresponding cosmooth projection
 valued measure on $\widehat{A}$ for $\pi_{1}$. Now an
 easy calculation show that the projection valued measure for the
 representation
 $\{a\longrightarrow\pi_{2}(h)\pi_{1}(a)\pi_{2}(h)^{-1}\}$
 of $A$ is $\{E\longrightarrow \pi_{2}(h)P_{E}\pi_{2}(h)^{-1}\}$,
 and that for the representation
 $\{a\longrightarrow\pi_{1}(h[a])\}$
 of $A$ is $\{E\longrightarrow P_{h[E]}\}$. In view of the
 uniqueness of the cosmooth projection valued measure which
 corresponds to a representation of $A$, we infer that
 $$\pi_{2}(h)P_{E}\pi_{2}(h)^{-1}=P_{h[E]}\;\;\;\eqno(2.11),$$
 so $(\pi_{2},P;V)$ is a cosmooth system of imprimitivity of
 $H$ based on $\widehat{A}$.

 Now let us state with $\pi_{1}$, $\pi_{2}$ and $P$ such that
$$\pi_{2}(h)P_{E}\pi_{2}(h)^{-1}=P_{h[E]},$$
then we gain (2.10). Define $\pi$ on $G$
$$\pi(ah)=\pi_{1}(a)\pi_{2}(h).$$
Then (2.10) is enough to secure the fact that $\pi$ is a
representation. Since the restriction of $\pi$ to $A$ and $H$ are
smooth and $G$ is equipped with the product topology, $\pi$ is
smooth, too.
\end{proof}

 The lemma just stated enables us to relate a smooth representation
 of $G$ to a cosmooth system of imprimitivity of $H$ based on
 $\widehat{A}$. The following lemma tells us that the relation is
 one to one in the sense of equivalence.

{\bf Lemma 2.} A smooth representation $\pi$ of $G$ on $V$ is
irreducible if and only if the corresponding cosmooth system of
imprimitivity for $H$ based on $\widehat{A}$ is irreducible. Two
smooth representations of $G$ are equivalent if and only if the
corresponding cosmooth systems of imprimitivity are equivalent.

\begin{proof}
For the first assertion we need only to prove that any subspace
$V_{1}$ of $V$ is invariant under $\pi_{1}=\pi|A$ if and only if
it is invariant under the corresponding cosmooth projection valued
measure $P$ based on $\widehat{A}$.

If $V_{1}$ is invariant under $P$, then by the definition (1.6),
$V_{1}$ is invariant under $\pi_{1}$. Now assume that $V_{1}$ is
invariant under $\pi_{1}$. Let $V^{\perp}_{1}$ be the subspace of
$\hat{V}$ (the space of linear functions on $V$) which is zero on
$V_{1}$. For any $v\in V_{1}, \widehat{v}\in V^{\perp}_{1}$,
$$<\pi_{1}(a)v, \widehat{v}>=0$$
Now
$$<\pi_{1}(a)v,\widehat{v}>=\int_{E(v)}x(a)<dP(x)v,\widehat{v}> \eqno(2.12),$$
and the basic fact in section 1.1 tell us that if $v\in V_{1}$,
then
$$<P(E)v,\widehat{v}>=<P(E\cap E(v))v,\widehat{v}>=0$$
for any open subset of $\widehat{A}$, and any $\widehat{v}$ in $
V^{\perp}_{1}$. Therefore $P(E)v\in  V_{1}$. In other word, $
V_{1}$ is invariant under $P$.

For the second assertion, let $\pi^{i}$ be smooth representations
of $G$ in $ V^{i}$, and let $(\pi_{2}^{i},P^{i}; V^{i})$ be the
corresponding cosmooth systems of imprimitivity, (i=1,2). Let $T$
be a isomorphism from $ V^{1}$ to $ V^{2}$. As in the proof of the
lemma 1, the cosmooth projection valued measure corresponding to
$T(\pi^{1}|A)T^{-1}$ is $TP^{1}T^{-1}$. Therefore by the
uniqueness stated in Theorem 1 shows that,
$\pi^{2}|A=T(\pi^{1}|A)T^{-1}$ if and only if
$P^{2}=TP^{1}T^{-1}$. From this the second assertion follows.
\end{proof}

Lemma 2 tells us that, to study smooth representations of $G$ is
equivalent to study the cosmooth systems of imprimitivity of
$(H,\widehat{A})$.

\subsection{Sheaf}
In this section, we use the concepts of presheaf and sheaf in the
sense of Bernstein and Zelevinsky.

We assume that $X$ is a totally disconnected locally compact space
and that $\mathcal{I}_{c}$ is the set of all compact open subsets
of $X$.

Let $C_{c}^{\infty}(X)$ be the sheaf of smooth complex valued
functions on $X$ with compact support. let $\mathcal{M}$ be a
sheaf of vector spaces over $X$ with base $\mathcal{I}_{c}$. Then
$\mathcal{M}$ is naturally a sheaf of module for
$C_{c}^{\infty}(X)$.

We call a $C_{c}^{\infty}(X)$-module $M$ {\bf{cosmooth}} if for
every $m \in M$, there exists a compact open subset $U$ of $X$
such that $1_{U}m=m$.

We have the following important proposition. For a proof, see
\cite{Bump}.

{\bf Proposition 1.} Let $M$ be a cosmooth
$C_{c}^{\infty}(X)$-module. We associate a presheaf $\mathcal{M}$
in the follow way. If $U\in \mathcal{I}_{c}$, let
$\mathcal{M}(U)=1_{U}\cdot M$. If $U\supseteq V,$ with $ U,V\in
\mathcal{I}_{c},$ we define a restriction map
$\rho_{U,V}:\mathcal{M}(U)\longrightarrow\mathcal{M}(V)$ by
$\rho_{U,V}(m)=1_{V}m$. Then $\mathcal{M}$ is a sheaf.

\subsection{Irreducible smooth representation}

Let $\pi$ be an irreducible admissible representation of $G$ and
let $(\pi_{2},P;V)$ be the corresponding cosmooth system of
imprimitivity of $H$ based on $\widehat{A}$ (see lemma 1). Due to
(1.5.1)-(1.5.4), $V$ becomes a cosmooth
$C_{c}^{\infty}(\widehat{A})$-module by setting the action of
$\chi_{E}$ on $V$ to be $P_{E}$. We can associate a sheaf
$\mathcal{V}$ to $V$ via proposition 1. $H$ has an action
$\Pi_{2}$ on $\mathcal{V}$ in the natural way, under which,
$\mathcal{V}(E)$ is mapped to $\mathcal{V}(h[E])$ and
$\mathcal{V}_{x}$ is mapped to $\mathcal{V}_{h[x]}$ by
$\Pi_{2}(h)$.

$\mathbf{Definition}$. Let $X$ be a $T_{1}$ $H$-space. $X$ is said
to be a {\bf{smooth $H$-space}} or in other word, $H$ acts
{\bf{smoothly}} on $X$, if for any two points $x_{1}$, $x_{2}$ in
$X$, either $x_{1}$ and $x_{2}$ lie in the same orbit of $H$ in
$X$, or there is $H$-invariant open subset of $X$ such that exact
one of $x_{1},x_{2}$ lies in it.

Now, we add a condition that $H$ acts smoothly on $\widehat{A}$.
Note that if $H$ is a compact group, then it acts always smoothly
on $\widehat{A}$.

{\bf Lemma 3.} Let $(\pi_{2},P;V)$, $(\pi_{2}^{1},P^{1};V^{1})$
and $(\pi_{2}^{2},P^{2};V_{2})$ be three irreducible cosmooth
systems of imprimitivity, and let $\mathcal{V}$, $\mathcal{V}^{1}$
and $\mathcal{V}^{2}$ be the sheaves associated to them
constructed by proposition 1. Then $supp(\mathcal{V})$ lies on
exact one orbit of $H$. If $supp(\mathcal{V}^{1})\neq
supp(\mathcal{V}^{2})$, then $(\pi_{2}^{1},P^{1};V^{1})$ and
$(\pi_{2}^{2},P^{2};V_{2})$  are two inequivalent cosmooth systems
of imprimitivity.

\begin{proof} We assert that $supp(\mathcal{V})$ lies on exact one
$H$-orbit. Otherwise, there are two orbits
$H\widehat{a_{1}},H\widehat{a_{2}}\subset supp(\mathcal{V})$, then
there is an $H$-invariant open subset $\tilde{E}$ such that exact
one of $H\widehat{a_{1}},H\widehat{a_{2}}$, say
$H\widehat{a_{1}}$, lies in $\tilde{E}$. Now, let
$\tilde{V}=P_{\tilde{E}}V$, then $\tilde{V}$ is a nontrivial
subspace of $V$, which is invariant under $(\pi_{2},P)$. This
contradicts the irreducibility of $(\pi_{2},P;V)$.

Let $(\pi^{1}_{2},P^{1};V^{1}), (\pi^{2}_{2},P^{2};V^{2})$ be two
irreducible cosmooth systems of imprimitivity. If
$supp(\mathcal{V}^{1}),supp(\mathcal{V}^{2})$ lie in two different
$H$-orbits $Hx_{1}, Hx_{2}$. Suppose $\tilde{E}$ is a
$H$-invariant open subset such that exact one of these two orbits
say $Hx_{1}$ lies in it, then $P^{1}_{\tilde{E}}\neq 0$, but
$P^{2}_{ \tilde{E} }=0$. Therefore $(\pi^{1}_{2},P^{1};V^{1}),
(\pi^{2}_{2},P^{2};V^{2})$ are two inequivalent cosmooth systems
of imprimitivity.
\end{proof}

Now let $(\pi_{2},P;V)$ be an irreducible cosmooth system of
imprimitivity, with $supp(\mathcal{V})$ lying in an orbit
$Hx_{0}$. Let $H_{0}$ be the stable subgroup of $x_{0}$ in $H$. It
is easy to see that $\mathcal{V}_{x_{0}}$, denoted by $V_{0}$, is
invariant under $H_{0}$. Let $\pi_{0}$ denote the action of
$H_{0}$ on $V_{0}$. By proposition 1, we can identity the sections
of $\mathcal{V}$ with the vectors in $V$. For every section $s$ of
$\mathcal{V}$, define a function on $H$ with value in $V_{0}$, by
$$\overrightarrow{F}_{s}(h)=\Pi_{2}(h)s(h^{-1}[x_{0}])\;\;\;\eqno\eqno(2.13).$$

Let $C^{\infty}_{c}(H/H_{0},\pi_{0},V_{0})$ denote the space of
locally constant functions $f$ with values in $V_{0}$ whose
support is compact $\mathrm{mod}\:H_{0}$, and satisfies
$$f(h_{0}h)=\pi_{0}(h_{0})f(h)\;\;\forall h_{0}\in H_{0}\;\forall h\in H\;\;\;\eqno\eqno(2.14).$$

{\bf Lemma 4.} $\overrightarrow{F_{s}}$ belongs to
$C^{\infty}_{c}(H/H_{0},\pi_{0},V_{0})$.

\begin{proof} It is easy to see that $\overrightarrow{F_{s}}$ satisfies
(2.14) and its support is compact $\mathrm{mod}\:H_{0}$. Let
$v_{s}$ be the vector in $V$ corresponding to $s$. Then there is a
compact open subgroup $H_{s}$ which fixes $v_{s}$. Note that
$$(\pi_{2}(h)s)(x)=\Pi_{2}(h)s(h^{-1}[x]).\;\;\;\eqno(2.15)$$

Therefore we have $\forall\;\; h\in H_{s},$
$$\Pi_{2}(h)s(h^{-1}[x])=s(x)\;\;\;\forall x\in \widehat{A}.$$
Especially, $\forall\;h\in H_{s}$
$$\overrightarrow{F}_{s}(gh)=\Pi_{2}(gh)s(h^{-1}g^{-1}[x_{0}])=\Pi_{2}(g)s(g^{-1}[x_{0}])=\overrightarrow{F}_{s}(g)$$
Thus $\overrightarrow{F}_{s}\in
C^{\infty}_{c}(H/H_{0},\pi_{0},V_{0} )$.
\end{proof}

Conversely, for a function $f\in
C^{\infty}_{c}(H/H_{0},\pi_{0},V_{0})$, we can define
$$\overleftarrow{F}_{f}(h[x_{0}])=\Pi_{2}(h)f(h^{-1})\;\;\;\eqno(2.16).$$
By (2.14), it is well defined.

{\bf Lemma 5.} $\overleftarrow{F}_{f}$ is a section of
$\mathcal{V}$.

\begin{proof}It is easy to see that it has compact support.

There is a compact open subgroup $H_{f}$ such that for each
$h_{f}\in H_{f}$, $f(hh_{f})=f(h)$. Then (2.16) tells us that
$$\Pi_{2}(h_{f})\overleftarrow{F}_{f}(x)=\overleftarrow{F}_{f}(h_{f}[x])\;\;\forall h_{f}\in H_{f},x\in Hx_{0}\;\;\;\eqno(2.17).$$

We are now to prove $\overleftarrow{F}_{f}$ is a section. Fix an
$x\in Hx_{0}$. We can select a section $s$ such that $s(x)=
\overleftarrow{F}_{f}(x) $. Let $H_{s}$ be a compact open subgroup
on $H$ such that $s$ is fixed by $H_{s}$. By (2.15), we gain
$$\Pi_{2}(h_{s})s(x)=s(h_{s}[x])\;\;\forall h_{s}\in H_{s},\; x\in Hx_{0} \eqno(2.18).$$
Comparing (2.17) and (2.18), we see that
$$s(h[x])= \overleftarrow{F}_{f}(h[x])\;\;\;\forall h\in H_{s}\cap H_{f}. $$
Thus $\overleftarrow{F}_{f}$ is really a section.
\end{proof}

It is obviously that $\overleftarrow{F}\overrightarrow{F}$ and
$\overrightarrow{F}\overleftarrow{F}$ are both identity, or
equivalently, $\overleftarrow{F}=\overrightarrow{F}^{-1}$.

Now $(\overrightarrow{F}\pi_{2}
\overleftarrow{F},\overrightarrow{F}P\overleftarrow{F};C_{c}^{\infty}(H/H_{0},\pi_{0},V_{0})
)$ is a cosmooth system of imprimitivity  that is equivalent to
$(\pi_{2},P;V)$. Write $(\bar{\pi}_{2},\bar{P};\bar{V})$ for
$(\overrightarrow{F}\pi_{2}
\overleftarrow{F},\overrightarrow{F}P\overleftarrow{F};C_{c}^{\infty}(H/H_{0},\pi_{0},V_{0})
)$. A direct calculation implies that: $\forall f \in \bar{V}$,
$$(\bar{\pi}_{2}(h_{1})f)(h)=f(hh_{1})\;\;\;\eqno(2.19.1),$$
$$\bar{P}_{E}f=\chi_{ \tilde{E}^{-1} }\cdot f\;\;\;\;\;\;\;\;\;\;\;\;\;\;\;\;\eqno(2.19.2),$$
where, $h_{1}\in H$, $\tilde{E}=\{h\in H: h[x_{0}]\in E\}$ and
$\tilde{E}^{-1}=\{h\in H; h^{-1}\in \tilde{E}\}$.

We see that $\overline{\pi}_{2}$ is just the compact induced
representation of $\pi_{0}$. Denote by $\overline{\pi}$, the
representation of $G$ corresponding to
$(\overline{\pi}_{2},\overline{P};\overline{V})$. A direct
calculation shows:
$$(\overline{\pi}(a)f)(h_{1})=<a, h_{1}^{-1}[x_{0}]>\cdot f(h_{1})\;\;\;\eqno(2.20.1),$$
$$(\bar{\pi}(h)f)(h_{1})=f(h_{1}h)\;\;\;\eqno(2.20.2),$$
where $a\in A$, $h\in H$ and $f\in \overline{V}$.

{\bf Lemma 6.}  The representation $\bar{\pi}$ is equivalent to
the representation $\pi$ mentioned  at the beginning of this
subsection.

This is just a consequence of Lemma 2.

Let $(\pi^{1}_{2},P^{1};V^{1})$ and $(\pi^{2}_{2},P^{2};V^{2})$ be
two irreducible cosmooth systems of imprimitivity, supported both
on $Hx_{0}$. Then it is easy to see that:
$$Hom((\pi^{1}_{2},P^{1};V^{1}),(\pi^{2}_{2},P^{2};V^{2}))\cong Hom((\Pi^{1},\mathcal{V}^{1}),(\Pi^{2},\mathcal{V}^{2}))$$
$$\cong Hom(\pi_{0}^{1},\pi_{0}^{2})\cong Hom(\bar{\pi}^{1},\bar{\pi}^{2}) \;\;\;\eqno(2.21)$$
Therefore $\bar{\pi}^{1}$ or $(\pi^{1},P^{1};V^{1})$ is
irreducible if and only if $\pi_{0}^{1}$ is irreducible. Moreover,
$(\pi^{1}_{2},P^{1};V^{1})$ and $(\pi^{2}_{2},P^{2};V^{2})$ are
two equivalent irreducible cosmooth systems of imprimitivity if
and only if $\pi^{1}_{0}$ and $\pi_{0}^{2}$ are two equivalent
irreducible smooth representations.

Write $\chi$ for $x_{0}$ now. For a representation $\pi_{0}$ of
$H_{\chi}$, let $\pi_{0}\cdot\chi$ be the representation of
$H_{\chi}\times_{t}A$:
$$(\pi_{0}\cdot\chi)(h\times a)=\chi(a)\pi_{0}(h)\;\;\forall\; h\in H_{\chi},\: a\in A.$$
It is easy to check that $\pi_{0}\cdot\chi$ is a representation.
Due to formulas (2.20.1) and (2.20.2), a simple calculate shows
that $\pi$ is equivalent to the compact induced representation
$\mathrm{Ind}^{G}_{H_{\chi}\times_{t} A}(\pi_{0}\cdot\chi)$ of
$G$.

Now, let $A$ be a locally compact and totally disconnected Abelian
group whose dual $\hat{A}$ having the same property. Let $H$ be a
locally compact and totally disconnected group with a continuous
action $t$ on $A$, and a dual action $t'$ on $\hat{A}$. Let $G$ be
$H\times_{t}A$.

We obtain the main result:

$\mathbf{Theorem}$ 2. For each orbit of $t'$, select out a point
$\chi$ on it. Every irreducible smooth representation $\pi_{0}$ of
$H_{\chi}$ gives an irreducible smooth representation
$\mathrm{Ind}^{G}_{H_{\chi}\times_{t} A}(\pi_{0}\cdot\chi)$ of
$G$. Every irreducible smooth representation of $G$ is equivalent
to one obtained in such a way. If furthermore $t'$ is smooth, the
representations obtained in such a way are not equivalent with
each other.


\begin{thebibliography}{9}

\bibitem[1]{Bernstein1} Bernstein, I.N., and Zelevinsky, A.V.: Induced representations
of reductive $p$-adic groups $\mathbf{I}$. Ann. Sci. Ecole Norm.,
Sup $4^{e}$ serie, t.10, 441-472 (1977)

\bibitem[2]{Bernstein2} Bernstein, I.N., and Zelevinsky, A.V.: Representations of the
group $GL(n,F)$ where $F$ is a local Non-Archimedean Field.
Russian Mathematical Surveys, no 3, 1-68 (1976)

\bibitem[3]{Bump} Bump, D.: Automorphic Forms and Representations,
Cambridge, Cambridge University Press, 1998

\bibitem[4]{Ja} Jacquet, H.: Generic representation. In
Non-Nommutative Harmonic Analysis, Lecture Notes in Mathematics,
587, 91-102 (1976)

\bibitem[5]{Mackey} Mackey, G.W.: Unitary Group representations in Physics,
Probability, and Number Theory. Addison-Wesley publishing company,
Inc., 1978, 1989

\bibitem[6]{Weil} Weil, A.: Sur la formule  de Siegel dans la th\'{e}orie des
groupes classiques. Acta Math 113, 1-87 (1965)

\bibitem[7]{V} Varadarajan, V.S.: Geometry of Quantum Theory, New
York, Springer-Verlag, 1985

\end{thebibliography}
\end{document}